\documentclass[11pt]{amsart}

\usepackage{amsmath}
\usepackage{amssymb}
\usepackage{graphicx}

\newtheorem{theorem}{Theorem}[section]
\newtheorem{corollary}[theorem]{Corollary}

\theoremstyle{definition}

\newtheorem{remark}[theorem]{Remark}

\newcommand{\bs}{\boldsymbol}
\newcommand{\cal}{\mathcal}
\newcommand{\fr}{\mathfrak}
\newcommand{\bb}{\mathbb}
\newcommand{\nl}{\lVert}
\newcommand{\nr}{\rVert} 

\numberwithin{equation}{section}

\title[On sharp Agmon-Miranda maximum principles]
{On sharp Agmon-Miranda maximum principles}

\author[G. Kresin]{Gershon Kresin}

\address[Gershon Kresin]{Department of Mathematics, Ariel University, Ariel 40700, Israel}
\email{{\tt kresin@ariel.ac.il}}

\author[V. Maz'ya]{Vladimir Maz'ya}

\address[Vladimir Maz'ya]{Department of Mathematical Sciences, University of Liverpool, M$\&$O Building, Liverpool, L69 3BX, UK}
\address{Department of Mathematics, Link\"oping University,SE-58183 Link\"oping, Sweden}
\address{RUDN University, 6 Miklukho-Maklay St., Moscow, 117198, Russia}   
\email{\tt vladimir.mazya@liu.se}



\begin{document}

\maketitle
\large
\centerline{\sl Dedicated to Shmuel Agmon with great respect}
\vspace{5mm}

\begin{abstract} 
In this  survey we formulate our results on different forms of maximum principles for linear elliptic equations and systems.
We start with necessary and sufficient conditions for validity of the classical maximum modulus principle for solutions of second order strongly elliptic systems. This  principle  holds under rather heavy restrictions on the  coefficients of the systems, for instance, it fails
for the Stokes and Lam\'e systems. Next, we turn to sharp constants in more general maximum principles due to S. Agmon and C. Miranda. 
We consider higher order elliptic equations, Stokes and  Lam\'e systems in a half-space as well as the system of planar deformed state in a half-plane.
\end{abstract}

\section{Introduction}\label{S_1}
Maximum principles are fundamental properties of partial differential operators, both linear and nonlinear. 
They have important applications to various facts in the theory of boundary value problems for these operators.

The present survey contains formulations of authors' results on the best constants in different forms of maximum principles for linear elliptic equations and systems.

Let $\Omega$ be a bounded domain in the Euclidean space ${\bb R}^n$.
Consider the uniformly elliptic equation
\begin{equation} \label{Eq_1.1}
\sum ^n_{j,k=1}a_{jk}(x)\frac{\partial ^2u}{\partial x_j\partial x_k}-\sum ^n _{j=1}a_j(x)
\frac{\partial u}{\partial x_j}-a_0(x)u=0\;,\;\;\;\;x\in\Omega,
\end{equation}
with bounded coefficients, positive-definite matrix $((a_{jk}(x)))$, and with $a_0(x)\geq 0$.
The following basic fact is called the classical maximum modulus principle.  

An arbitrary solution $u\in {\rm C}^2({\Omega})\cap {\rm C}(\overline\Omega)$ of equation (\ref{Eq_1.1}) 
satisfies
\begin{equation} \label{Eq_1.2}
\max _{\overline \Omega}\; | u |\leq \max _{\partial \Omega}\;| u | \;.
\end{equation}

This principle was obtained for different classes of solutions and under various assumptions about the coefficients.

Henceforth by $|\cdot|$ we denote the Euclidean length of a vector. 
If, instead of equation (\ref{Eq_1.1}), we consider a homogeneous second order elliptic system with vector-valued solutions
$\bs u\in [{\rm C}^2({\Omega})]^m\cap [{\rm C}(\overline\Omega)]^m$ subject to the inequality 
\begin{equation} \label{Eq_1.2A}
\max _{\overline \Omega}\; | \bs u |\leq \max _{\partial \Omega}\;| \bs u | \;,
\end{equation}
then we say that the classical maximum modulus principle holds for the system in question. 

In section \ref{S_2} we state criteria for validity of the classical maximum modulus principle 
for second order strongly elliptic systems following our papers \cite{KM3,MK}.
Note that this principle holds for the systems under rather heavy conditions on the coefficients.
In particular, Polya's paper \cite{PO} contains an example showing that the classical maximum modulus principle fails for
the displacement vector satisfying the Lam\'e system. The same is true for the velocity vector subject to the Stokes system.

The not so sharp as the classical modulus principle but incomparably more general properties of the same nature, which hold for elliptic equations of arbitrary order in smooth domains, were discovered by S. Agmon and C. Miranda. 
Similar results for elliptic systems were obtained by S. Agmon, A. Douglis, L. Nirenberg.

For solutions of a homogeneous elliptic equation of order $2\ell$, Agmon-Miranda principle is the estimate
\begin{equation} \label{Eq_1.3}
\max _{\overline \Omega}\; | \nabla _{\ell-1}u | \leq c(\Omega)
\max _{\partial \Omega}\; | \nabla_{\ell-1} u | .
\end{equation}
A weaker variant of the Agmon-Miranda maximum principle runs as follows:
\begin{equation} \label{Eq_1.4}
\max _{\overline\Omega} |\nabla _{\ell-1} u|\leq K
\max _{\partial\Omega}|\nabla _{\ell-1} u|+C\nl u \nr_{L^1(\Omega)}.
\end{equation}
An estimate similar to (\ref{Eq_1.3}) for the biharmonic equation with two variables was proved earlier by C. Miranda \cite{MIR1}.
For strongly elliptic equations with real coefficients in the two-dimensional case, a result of type (\ref{Eq_1.3}) was obtained by C. Miranda
in \cite{MIR2} with the help of Agmon's result in \cite{AG0}. 
In the general case of  equation with complex coefficients with any number of independent variables, (\ref{Eq_1.3}) 
was established in Agmon \cite{AG}. 

The inequality 
\begin{equation} \label{Eq_1.5A}
\max _{\overline \Omega}\; |{\bf u}| \leq k(\Omega)
\max _{\partial \Omega}\;|{\bf u}|,\;\;\;\;\;\;\;\;\;\;\;
\Omega\subset {\bb R}^{n} 
\end{equation}
was established by Fichera \cite{FI} for solutions of the Lam\'e system and is called the Fichera's maximum principle in elastostatics.
The construction for $k(\Omega)$ proposed by Fichera, gives some information about
dependence of $k(\Omega)$ on the geometry of $\Omega$, and elastic constants $\lambda $ and $\mu $.
More general estimates for solutions of elliptic systems were obtained in Agmon, Douglis and Nirenberg \cite{ADN2}. 

In section \ref{S_3}, following our paper \cite{KM2}, we describe a sharp constant $K$ in inequality (\ref{Eq_1.4})  for solutions to higher order elliptic equation with constant complex coefficients. Besides, we give an explicit formula for the sharp constant $K$ in Miranda's inequality 
$$
\sup _{\overline {{\bb R}^n _{+}}} |\nabla u|\leq K
\sup _{\partial {\bb R}^n _{+}} |\nabla  u|
$$
for a biharmonic function $u$ in the half-space ${\bb R}^n _{+}=\{ x=(x_1,\dots, x_n)\in{\bb R}^n: x_n>0 \}$. 

In section \ref{S_4} we collect the results obtained in \cite{MK}.
Here we state explicit formulas for the sharp constant ${\mathcal K}({\bb R}^n _{+})$ in inequality 
\begin{equation} \label{Eq_1.6A}
|\bs u(x)|\leq {\cal K}({\bb R}^n _{+})\sup \{|\bs u(x')|: x' \in \partial {\bb R}^n _{+} \},
\end{equation}
for bounded solutions of the Lam\'e and Stokes systems in a half-space. 
Besides, we give sharp majorant for two components  of the stress tensor of the planar deformed state in a half-plane. 

\smallskip
The above mentioned authors' results with detailed proofs have been collected in monograph \cite{KM}.  

\section{Classical maximum modulus principle for solutions to second order strongly elliptic systems}\label{S_2}

The first articles aiming at the study and applications of the classical maximum modulus principle for solutions of elliptic second order systems concerned systems with scalar coefficients in the first and second derivatives of the unknown vector-valued function (see Bitsadze \cite{BI,BITS}, 
Pini \cite{PI}, Szeptycki \cite{SZE}). These systems are weakly coupled: a system of partial differential equations is called weakly coupled if there are no derivatives in the coupling terms.

Sufficient conditions for validity of the maximum modulus principle, its modifications
and generalizations for non-weakly coupled systems were given by
Hile and Protter \cite{HP}, C. Miranda \cite{MIR3},
Rus \cite{Rus}, Stys \cite{STY1}, Sabitov \cite{SAB2}, Wasowski \cite{WAS}. In particular, C. Miranda \cite{MIR3}
considered elliptic second order systems with a scalar principal part and with arbitrary coefficients 
in derivatives of order less than two. He found an algebraic inequality sufficient 
for the classical maximum modulus principle
(conditions in Remarks \ref{R1} and \ref{R3} with  the strict inequality sign). 
A survey of maximum principles for elliptic equations and systems with a scalar 
principal part is given by Protter \cite{PR}.

An algebraic necessary and sufficient condition for validity of the maximum
principle for the product $\alpha(x)|\bs u|$, where $\alpha$ is a certain function and
$\bs u$ is a solution of the elliptic system with analytic coefficients, 
is due to Hong \cite{Hong}. 

There is a number of results on the componentwise maximum principle for weakly 
coupled elliptic systems, in particular, on non-negativity of the components of a solution
(see de Figueiredo and Mitidieri \cite{FiMi},  Lenhart and Schaefer \cite{LeS},
L\'opez-G\'omez and Molina-Meyer \cite{LoMo},  Mitidieri and Sweers \cite{MiSw}, Sirakov \cite{SirB} 
and bibliography there). Various maximum principles  for weakly coupled systems
are discussed in the book by Protter and Weinberger \cite{PW}.
Necessary and sufficient conditions for the componentwise and for the so-called "stochastic" 
extremum principles for solutions of elliptic systems of the second order are given by 
Kamynin and Khimchenko in \cite{KamKh1} and \cite{KamKh2}, respectively.

\smallskip
In this section we describe criteria for validity of the classical maximum modulus principle for solutions of the strongly elliptic system
\begin{equation} \label{(3.1)}
\sum ^n_{j,k=1}{\cal A}_{jk}(x)\frac{\partial ^2\bs u}{\partial x_j\partial x_k}-\sum ^n _{j=1}{\cal A}_j(x)\frac{\partial \bs u}{\partial x_j}-{\cal A}_0(x)\bs u=\bs 0 
\end{equation}
with real or complex coefficients. Here ${\cal A}_{jk}, {\cal A}_{j},
{\cal A}_{0}$ are $(m\times m)$-matrix-valued functions and $\bs u$ is a
$m$-component vector-valued function. Without loss of generality we assume that ${\cal A}_{jk}={\cal A}_{kj}$.

\subsection{Model systems}

We begin with the simple case of the homogeneous operator with the constant coefficients.

\subsubsection{The case of real coefficients}

We introduce the operator
\[
{\fr A}_{0}(D_x)=\sum ^n_{j,k=1}{\cal
A}_{jk}\frac{\partial ^2}{\partial x_j\partial x_k},
\]
where $D_x=(\partial /\partial x_1,\dots,\partial /\partial x_n)$ and 
${\cal A}_{jk}={\cal A}_{kj}$ are constant real $(m\times m)$-matrices. Assume that the
operator ${\fr A}_{0}$ is strongly elliptic, i.e. that for all $\bs\zeta=(\zeta_1,
\dots,\zeta_m)\in {\bb R}^m$ and $\bs\sigma=(\sigma _1,\dots,\sigma _n)\in
{\bb R}^n$, with $\bs\zeta,\;\bs\sigma \neq \bs 0$, we have the inequality
\[
\biggl (\;\sum^n _{j,k=1}{\cal A}_{jk}\sigma _{j}
\sigma _{k}\bs\zeta,\; \bs\zeta 
\biggr )>0.
\]

Let $\Omega$ be a domain in ${\bb R}^n$ with boundary $\partial \Omega$
and closure $\overline \Omega$. Let $[{\rm C}_{\rm b}({\overline \Omega})]^m$ denote the space
of bounded $m$-component vector-valued functions which are continuous in
$\overline \Omega$. The norm on $[{\rm C}_{\rm b}({\overline \Omega})]^m$ is $\nl\bs u\nr=\sup \big
\lbrace |\bs u(x)|:\;x\in {\overline \Omega}\big \rbrace$. The notation
$[{\rm C}_{\rm b}({\partial \Omega})]^m$ has a similar meaning. 
By $[{\rm C}^{2}(\Omega)]^m$ we denote the space of $m$-component vector-valued 
functions with continuous derivatives up to the second order in $\Omega$. 
We omit the upper index $m$ in notation of spaces in the case $m=1$.

Let 
\begin{equation} \label{EELMSR_0.01}
{\cal K}(\Omega)=\sup \frac{\nl\bs u\nr_{[{\rm C}_{\rm b}(\overline \Omega)]^m}}
{\nl\bs u\nr_{[{\rm C}_{\rm b}(\partial \Omega )]^m}},
\end{equation}
where the supremum is taken over all vector-valued functions in the class
$[{\rm C}_{\rm b}(\overline \Omega)]^m \cap [{\rm C}^{2}(\Omega)]^m$ satisfying the system
${\fr A}_{0}(D_x)\bs u=\bs 0$.

Clearly, ${\cal K}(\Omega)$ is the best constant in the inequality
\begin{equation} \label{Eq_Main}
|\bs u(x)| \leq {\cal K}(\Omega)\sup\;\lbrace |\bs u(y)| : y\in \partial \Omega \rbrace,
\end{equation}
where $x\in \Omega$ and $\bs u$ is a solution of the system
${\fr A}_{0}(D_x)\bs u=0$
in the class $[{\rm C}_{\rm b}(\overline \Omega)]^m \cap [{\rm C}^{2}(\Omega)]^m$.

If ${\cal K}(\Omega)=1$, then the classical maximum modulus principle holds for
the system ${\fr A}_{0}(D_x)\bs u=\bs 0$.

According to Agmon, Douglis and Nirenberg \cite{ADN2},  Lopatinski\v{\i} \cite{LO1}, Shapiro \cite{SHA}, Solonnikov \cite{SOL} 
there exists a bounded solution of the Dirichlet problem
\begin{equation} \label{EELMSR_0.02}
{\fr A}_{0}(D_x)\bs u=\bs 0\; \hbox {in}\; {\bb R}^n_{+},\;\; \bs u=\bs f
\; \hbox {on}\; \partial {\bb R}^n _{+},
\end{equation}
with $\bs f\in [{\rm C}_{\rm b}(\partial {\bb R}^n _{+})]^m$, such that $\bs u$ is continuous
up to $\partial {\bb R}^n _{+}$, and can be represented in the form
\begin{equation} \label{EELMSR_0.03}
\bs u(x)=\int _{\partial {\bb R}^n_{+}} {\cal M}  \left (\frac{y-x}{|y-x|}\right )
\frac{x_n}{|y-x|^n}\bs f(y')dy'.
\end{equation}
Here $y=(y', 0),\; y'=(y_1,\dots,y_{n-1}),$ and ${\cal M}$ is a continuous
$(m\times m)$-matrix-valued function on the closure of the hemisphere
${\bb S}^{n-1} _{-}=\big \lbrace x\in {\bb R}^n :\;|x|=1, x_{n} < 0 \big \rbrace$
such that  the integral
\[
\int _{{\bb S}^{n-1} _ {-}}{\cal M}(\sigma ) d \sigma
\]
is the identity matrix.

The uniqueness of the solution of the Dirichlet problem (\ref{EELMSR_0.02}) in the
class $[{\rm C}_{\rm b}(\overline {{\bb R}^n _{+}})]^m \cap [{\rm C}^{2}({\bb R}^n _{+})]^m$ 
can be derived by means of standard arguments from (\ref{EELMSR_0.03}) and by local estimates
of the derivatives of solutions to elliptic systems (see Agmon, Douglis and Nirenberg \cite{ADN2}, Solonnikov \cite{SOL}).

The  assertion below contains a representation of the sharp constant ${\cal K}({\bb R}^n _{+})$ in the pointwise estimate (\ref{Eq_Main}) 
for the case $\Omega={\bb R}^n _{+}$. 
\setcounter{theorem}{0}
\begin{theorem} \label{TELMSR_0.01}
The formula
\begin{equation} \label{EELMSR_0.05}
{\cal K}({\bb R}^n _{+})=\sup_{|\bs{z}|=1}\int _{{\bb S}^{n-1} _{-}}|{\cal M}^*(\sigma)\bs{z}|d\sigma
\end{equation}
is valid, where asterisk denotes passage to the transposed matrix and $\bs z\in {\bb R}^m$.
\end{theorem}

The next statement gives a criterion for validity of the classical modulus principle for the system 
\begin{equation} \label{MSRC} 
\sum ^n_{j,k=1}
{\cal A}_{jk}\frac{\partial ^2\bs u}{\partial x_j\partial x_k}=\bs 0  
\end{equation}
in ${\bb R}^n_+$ with the real constant coefficients.
\begin{theorem} \label{TELMSR_0.02}
The equality ${\cal K}({\bb R}^n _{+})=1$ is satisfied if and only if
\begin{equation} \label{EELMSR_0.07}
{\fr A}_{0}(D_x)={\cal A}\sum ^n_{j,k=1}
a_{jk}\frac{\partial ^2}{\partial x_j\partial x_k},
\end{equation}
where ${\cal A}$ and $((a_{jk}))$ are positive-definite constant matrices of orders $m$ and $n$, respectively.
\end{theorem}

The theorem below, together with Theorem \ref{TELMSR_0.02} form a necessary condition for the 
validity of the classical maximum modulus principle for system (\ref{MSRC}). 
\begin{theorem} \label{TELMSR_0.03}
Let $\Omega$ be a domain in ${\bb R}^n$ with
compact closure and ${\rm C}^1$-boundary. Then
\[
{\cal K}(\Omega )\geq\sup\lbrace {\cal K}({\bb R}^n_{+}(\bs\nu)) :
\bs\nu\in {\bb S}^{n-1}\rbrace,
\]
where ${\bb R}^n_{+}(\bs\nu)$ is the half-space with inward normal $\bs\nu$
and ${\bb S}^{n-1}=\{x\in {\bb R}^n: |x|=1 \}.$ 
\end{theorem}

Now, we formulate a necessary and sufficient condition ensuring the classical modulus principle
for system (\ref{MSRC}) in a bounded domain with a smooth boundary.
\begin{theorem} \label{TELMSR_0.04}
Let $\Omega $ be a domain in ${\bb R}^n$
with compact closure and ${\rm C}^1$-boundary. The equality ${\cal K}(\Omega)=1$
holds if and only if the operator ${\fr A}_{0}(D_x)$ is defined by $(\ref{EELMSR_0.07})$.
\end{theorem}

\subsubsection{The case of complex coefficients}

We introduce the operator
\[
{\fr C}_{0}(D_x)=\sum ^n_{j,k=1}{\cal C}_{jk}
\frac{\partial ^2}{\partial x_j\partial x_k},
\]
where ${\cal C}_{jk}={\cal C}_{kj}$ are constant complex $(m\times m)$-matrices. Assume that the
operator ${\fr C}_{0}$ is strongly elliptic, that is
\[
\Re \left ( \sum^n _{j,k=1}{\cal C}_{jk}\sigma _{j}
\sigma _{k}\bs\zeta,\; \bs\zeta \right ) > 0
\]
for all $\bs\zeta=(\zeta_1,\dots,\zeta_m)\in {\bb C}^m$ and $\bs\sigma=(\sigma _1,\dots,
\sigma _n)\in {\bb R}^n$, with $\bs\zeta,\;\bs\sigma \neq \bs 0$.
Here and henceforth ${\bb C}^m$ is a complex linear $m$-dimensional space
with the elements $\bs a+i\bs b$, where $\bs a, \bs b\in {\bb R}^m$. The inner product
in ${\bb C}^m$ is $\left (\bs c, \bs d \right )=c_1{\overline d_1}+\dots+
c_m{\overline d_m},\;\bs c=(c_1,\dots,c_m),\;\bs d=(d_1,\dots,d_m)$. The
length of the vector $\bs d$ in ${\bb C}^m$ is $|\bs d|=(\bs d, \bs d )^{1/2}$.

\medskip
Let ${\cal R}_{jk}$ and ${\cal H}_{jk}$ be constant real $(m\times m)$-matrices
such that ${\cal C}_{jk}={\cal R}_{jk}+i{\cal H}_{jk}$, $i=\sqrt{-1}$. We define the operators
\[
{\fr R}_{0}(D_x)=\sum ^n_{j,k=1}{\cal
R}_{jk}\frac{\partial ^2}{\partial x_j\partial x_k},\;\;\;\;
{\fr H}_{0}(D_x)=\sum ^n_{j,k=1}{\cal
H}_{jk}\frac{\partial ^2}{\partial x_j\partial x_k}.
\]
Separating the real and imaginary parts of the system ${\fr C}_{0}(D_x)\bs u=0$, 
where $\bs u=\bs v+i\bs w$, we get a system with real coefficients
\[
{\fr R}_{0}(D_x)\bs v-{\fr H}_{0}(D_x)\bs w=\bs 0,\;\;\;\;
{\fr H}_{0}(D_x)\bs v+{\fr R}_{0}(D_x)\bs w=\bs 0,
\]
which, like the original system, is strongly elliptic.

We introduce the matrix operator
\[
{\fr K}_{0}(D_x)=\begin{pmatrix}{\fr R}_{0}(D_x)& - {\fr H}_{0}(D_x)\\
  {\fr H}_{0}(D_x) & {\fr R}_{0}(D_x) \end{pmatrix}.
\]

Let $[{\bf C}_{\bf b}(\Omega)]^m$ be the space of $m$-component complex
vector-valued functions $\bs u=\bs v+i\bs w$ which are bounded and continuous on $\Omega
\subset {\bb R}^n$. The norm on $[{\bf C}_{\bf b}(\Omega)]^m$ is $\nl\bs u\nr=\sup \big \lbrace
(|\bs v(x)|^{2}+|\bs w(x)|^{2})^{1/2}\;:\;x\in \Omega \big \rbrace$.
The notation $[{\bf C}_{\bf b}({\partial \Omega})]^m$ has a similar meaning.
By $[{\bf C}^{2}(\Omega)]^m$ we denote the space of $m$-component complex
vector-valued functions with continuous derivatives up to the second
order in $\Omega$.

By analogy with the definition (\ref{EELMSR_0.01}) of ${\cal K}(\Omega)$, let
\[
{\cal K}'(\Omega)=\sup \frac{\nl\bs u\nr_{[{\bf C}_{\bf b}(\overline \Omega)]^m}}{\nl\bs u\nr_{[{\bf C}_{\bf b}(\partial \Omega)}]^m},
\]
where the supremum is extended over all vector-valued functions in the class
$[{\bf C}^{2}(\Omega)]^m\cap [{\bf C}_{\bf b}(\overline \Omega)]^m$ subject to the system
${\fr C}_{0}(D_x)\bs u=0$ in $\Omega$.

It is clear that the constant ${\cal K}'(\Omega)$ for the system
${\fr C}_{0}(D_x)\bs u=0$ with complex coefficients coincides
with the constant ${\cal K}(\Omega)$ for the system ${\fr K}_{0}(D_x)\big
\lbrace \bs v, \bs w \big \rbrace=0$ with real coefficients
if we replace $m$ by $2m$, ${\fr A}_{0}(D_x)$ by
${\fr K}_{0}(D_x)$, and $\bs u$ by $\big \lbrace \bs v, \bs w \big \rbrace$
in definition (\ref{EELMSR_0.01}). Therefore, all assertions about ${\cal K}'(\Omega)$
are direct consequences of the analogous assertions about ${\cal K}(\Omega)$.
Using this fact, for solutions of the system
\begin{equation} \label{MSCC} 
\sum ^n_{j,k=1}
{\cal C}_{jk}\frac{\partial ^2\bs u}{\partial x_j\partial x_k}=\bs 0  
\end{equation}
with constant complex coefficients, we obtain the following four theorems.
 
\setcounter{theorem}{4}
\begin{theorem} \label{TELMSC_0.01}
The formula
\[
{\cal K}'({\bb R}^n _{+})=\sup_{|\bs z|=1}\int _{{\bb S}^{n-1} _{-}}|{\cal U}^*
(\sigma)\bs z|d\sigma
\]
is valid, where ${\cal U}$ is the $(2m\times 2m)$-matrix-valued function on ${\bb S}^{n-1} _{-}$
appearing in the integral representation for a
solution of the Dirichlet problem in ${\bb R}^n _{+}$ for the system
${\fr K}_{0}(D_x)\big\lbrace \bs v, \bs w \big \rbrace=0$
$($analogous to representation $(\ref{EELMSR_0.03}))$ and $\bs z\in {\bb R}^{2m}$.
\end{theorem}
\begin{theorem} \label{TELMSC_0.02}
The equality ${\cal K}'({\bb R}^n _{+})=1$ holds if and only if
\begin{equation} \label{EELMSC_0.01}
{\fr C}_{0}(D_x)={\cal C}\sum ^n_{j,k=1}a_{jk}
\frac{\partial ^2}{\partial x_j\partial x_k},
\end{equation}
where ${\cal C}$ is a constant complex-valued
$(m\times m)$-matrix such that $\Re ({\cal C}\bs\zeta, \bs\zeta ) > 0$ 
for all $\bs\zeta\in {\bb C}^m,\;\bs\zeta\neq \bs 0$, and
$((a_{jk}))$ is a real positive-definite $(n\times n)$-matrix.
\end{theorem}

\begin{theorem} \label{TELMSC_0.03}
Let $\Omega$ be a domain in ${\bb R}^n$ with compact closure and 
${\rm C}^1$-boundary. Then
\[
{\cal K}'(\Omega)\geq\sup\lbrace {\cal K}'({\bb R}^n_{+}(\bs\nu)) :
\bs\nu \in {\bb S}^{n-1}\rbrace,
\]
where ${\bb R}^n_{+}(\bs\nu)$ is a half-space with inward normal $\bs\nu$.
\end{theorem}

A necessary and sufficient condition for validity of the classical modulus principle
for system  (\ref{MSCC}) in a bounded domain runs as follows. 
\begin{theorem} \label{TELMSC_0.04} Let $\Omega$ be a domain in
${\bb R}^n$ with compact closure and ${\rm C}^1$-boundary.
The equality ${\cal K}'(\Omega)=1$ holds if and only if
the operator ${\fr C}_{0}(D_x)$ has the form $(\ref{EELMSC_0.01})$.
\end{theorem}

\subsection{Systems with lower order terms and variable coefficients}

Now, we turn to linear elliptic system (\ref{(3.1)}) of the general form.

\subsubsection{The case of real coefficients} \label{SELMGR_1.01}

Let $\Omega$ be a domain in ${\bb R}^n$ with compact closure $\overline\Omega$
and with boundary $\partial\Omega$ of the class ${\rm C}^{2, \alpha}, 0 <\alpha \leq 1.$ 
The space of $(m\times m)$-matrix-valued functions whose
elements have continuous derivatives up to order $k$ and satisfy the H\"older
condition with exponent $\alpha, 0 <\alpha \leq 1$, on $\overline\Omega$  is denoted
by $[{\rm C}^{k, \alpha }(\overline\Omega)]^{m\times m}$. 

We introduce the operator
\[
{\fr A}(x, D_x)=\sum ^n_{j,k=1}{\cal A}_{jk}(x)
\frac{\partial ^2}{\partial x_j\partial x_k}-\sum ^n _{j=1}{\cal A}_j(x)
\frac{\partial}{\partial x_j}-{\cal A}_0(x),
\]
where ${\cal A}_{jk}={\cal A}_{kj}, {\cal A}_{j}, {\cal A}_{0}$ are real 
$(m\times m)$-matrix-valued functions in the spaces 
\[
[{\rm C}^{2, \alpha }(\overline\Omega)]^{m\times m},\;\;\;\;\;\;
[{\rm C}^{1, \alpha }(\overline\Omega)]^{m\times m},\;\;\;\;\;\; 
[{\rm C}^{\alpha }(\overline\Omega)]^{m\times m},
\]
respectively. If the coefficients of the operator ${\fr A}(x,D_x)$ 
do not depend on $x$ we use the notation ${\fr A}(D_x)$. Let the 
principal homogeneous part of the operator ${\fr A}(x, D_x)$ be denoted by 
${\fr A}_{0}(x, D_x)$.

We assume that ${\fr A}(x, D_x)$ is strongly elliptic
in $\overline\Omega$, which means that for all $x\in\overline\Omega,\;\bs \zeta=(\zeta_1,
\dots,\zeta_m)\in {\bb R}^m,\;\bs \sigma=(\sigma _1,\dots,\sigma _n)\in
{\bb R}^n$, with $\bs \zeta,\;\bs \sigma \neq \bs 0$, the inequality
\begin{equation} \label{EELMGR_1.01}
\Biggl (\sum^n _{j,k=1}{\cal A}_{jk}(x)\sigma _{j}
\sigma _{k}\bs\zeta,\; \bs\zeta \Biggr ) >0
\end{equation}
is satisfied.

The next assertion gives necessary and sufficient conditions
for validity of the classical maximum modulus principle for 
system (\ref{(3.1)}) in any subdomain $\omega$ of a bounded domain $\Omega$ with smooth boundary.

\begin{theorem} \label{TELMGR_1.07}
The classical maximum modulus principle
\begin{equation} \label{CMP}
\nl\bs u\nr_{[{\rm C}(\overline\omega )]^m} \leq  \nl\bs u|_{\partial\omega}
\nr_{[{\rm C}(\partial \omega )]^m},
\end{equation}
holds for solutions of the system ${\fr A}(x, D_x)\bs u=\bs 0$
in an arbitrary domain $\omega\subset\Omega$  with boundary
from the class  ${\rm C}^{2,\alpha}$ if and only if:

{\rm (i)} for all $x\in \overline\Omega$ the equalities hold
\[
{\cal A}_{jk}(x)={\cal A}(x)a_{jk}(x),\;\; 1\leq j,\; k\leq n,
\]
where ${\cal A}$ and $((a_{jk}))$ are real positive-definite matrices
in $\overline\Omega$ of orders $m$ and $n$, respectively;

{\rm (ii)} for all $x\in \Omega$  and any $\bs \xi_j,\;\bs\zeta \in {\bb R}^m,
\; j=1,\dots ,n,$  with $(\bs\xi _j, \bs\zeta)=0$ the inequality
\[
\sum ^n_{j,k=1}a_{jk}(x)(\bs\xi _j, \bs\xi _k)+\sum ^n_{j=1}({\cal A}^{-1}(x)
{\cal A}_j(x)\bs\xi_j, \bs\zeta )+({\cal A}^{-1}(x){\cal A}_0(x)\bs\zeta , \bs\zeta ) \geq 0
\]
is valid.
\end{theorem}

The next assertion is a consequence of Theorem \ref{TELMGR_1.07}.

\setcounter{theorem}{0}
\begin{corollary} \label{CELMGR_1.04}
The classical maximum modulus principle\index{classical maximum modulus principle}
\[
\nl\bs u\nr_{[{\rm C}(\overline \omega)]^m} \leq
\nl\bs u |_{\partial \omega }\nr_{[{\rm C}(\partial  \omega)]^m}
\]
holds for solutions of the system
\[
\sum_{j,k=1}^{n} {\cal A}_{jk}(x)\frac{\partial ^2 \bs u}{ \partial x_j \partial x_k}-
\sum_{j=1}^{n} {\cal A}_j(x)\frac{\partial \bs u}{\partial x_j}=\bs 0
\]
in an arbitrary domain $\omega \subset \Omega $ with boundary from the class ${\rm C}^{2, \alpha }$
if and only if:
\[
{\cal A}_{jk}(x)={\cal A}(x)a_{jk}(x),\;\;\;{\cal A}_{j}(x)={\cal A}(x)a_{j}(x),\;\; 1\leq j,\; k\leq n.
\]
Here ${\cal A}$ and $((a_{jk}))$ are positive-definite matrix-valued functions in 
$\overline \Omega $ of orders $m$ and $n$, respectively, and $a_j$ are scalar 
functions with the same smoothness as ${\cal A}_j$.
\end{corollary}

\setcounter{theorem}{0}   
\begin{remark}\label{R1}  Condition (ii) of Theorem \ref{TELMGR_1.07} can be replaced by the following:

{\it for all $x\in \Omega$ and for any
$\bs\zeta \in {\bb R}^m\backslash\lbrace \bs 0 \rbrace$ the inequality holds:
\begin{eqnarray*}
& &|\bs\zeta|^{-2}\sum ^n_{j,k=1}b_{jk}(x)({\cal A}^{-1}(x){\cal A}_j(x)
\bs\zeta , \bs\zeta)({\cal A}^{-1}(x){\cal A}_k(x)\bs\zeta , \bs\zeta ) \\
& &- \sum ^n_{j,k=1}b_{jk}(x)({\cal A}_j^*(x)({\cal A}^*(x))^{-1}\bs\zeta,\;{\cal A}_k^*(x)
({\cal A}^*(x))^{-1}\bs\zeta )\\
& &\\
& &+4({\cal A}^{-1}(x){\cal A}_0(x)\bs\zeta, \bs\zeta )\geq 0. 
\end{eqnarray*}
Here $((b_{ij}))$ is the inverse matrix of $((a_{ij}))$
and * means passage to the transposed matrix.}
\end{remark}

\begin{remark}\label{R2}
In \cite{KM3} we showed by an example  that the
possibility to represent the principal part of the system ${\fr A}(x, D_x)\bs u=\bs 0$ in the form
$$
{\cal A}(x)\sum ^n_{j,k=1} a_{jk}(x)\frac{\partial ^2
\bs u}{\partial x_j\partial x_k}
$$
everywhere in $\overline\Omega$ is not necessary for validity of the classical maximum 
modulus principle 
$$
\nl\bs u\nr_{[{\rm C}(\overline\Omega )]^m} \leq  \nl\bs u|_{\partial\Omega}
\nr_{[{\rm C}(\partial \Omega )]^m},
$$
where $\bs u$ is a solution of the system ${\fr A}(x, D_x)\bs u=\bs 0$ in $\Omega$ which belongs to
$[{\rm C}^2(\Omega)]^m\cap [{\rm C}(\overline {\Omega})]^m$.
\end{remark}

\subsubsection{The case of complex coefficients}

In this section we extend basic results of subsection \ref{SELMGR_1.01} to system
(\ref{(3.1)}) with complex coefficients with solutions $\bs u=\bs v+i\bs w$,
where $\bs v$ and $\bs w$ are $m$-component vector-valued functions with real-valued
components.
Here, similarly to subsection \ref{SELMGR_1.01}, we assume that
$\Omega$ is a domain in ${\bb R}^n$ with compact closure $\overline\Omega$
and with boundary $\partial\Omega$ in the class ${\rm C}^{2, \alpha}, 0 <\alpha \leq 1$. 

For the spaces of matrix-valued functions with complex
components we retain the same notation as in the case of real components
but use bold letters.

\smallskip
We introduce the  operator
\[
{\fr C}(x, D_x)=\sum ^n_{j,k=1}{\cal C}_{jk}(x)
\frac{\partial ^2}{\partial x_j\partial x_k}-\sum ^n _{j=1}{\cal C}_j(x)
\frac{\partial}{\partial x_j}-{\cal C}_0(x),
\]
where ${\cal C}_{jk}={\cal C}_{kj}, {\cal C}_{j}, {\cal C}_{0}$ are complex
$(m\times m)$-matrix-valued functions in the spaces
\[
[{\bf C}^{2, \alpha }(\overline\Omega)]^{m\times m},\;\;\;\;\;
[{\bf C}^{1, \alpha }(\overline\Omega)]^{m\times m},\;\;\;\;\;
[{\bf C}^{\alpha }(\overline\Omega)]^{m\times m},
\]
respectively.
Suppose that the operator ${\fr C}(x, D_x)$ is strongly
elliptic in $\overline\Omega$, that is for all $x\in\overline\Omega,\;
\bs\zeta=(\zeta_1,\dots,\zeta_m)\in {\bb C}^m,\;\bs\sigma=(\sigma _1,\dots,
\sigma _n)\in {\bb R}^n$, with $\bs\zeta,\;\bs\sigma \neq \bs 0$, the inequality
\[
\Re \Biggl ( \sum^n _{j,k=1}{\cal C}_{jk}(x)\sigma _{j}
\sigma _{k}\bs\zeta,\;\bs\zeta \Biggr )\; >0
\]
holds.

Let ${\cal R}_{jk}, {\cal H}_{jk}, {\cal R}_{j}, {\cal H}_{j}, {\cal R}_{0},
{\cal H}_{0}$ be real $(m\times m)$-matrix-valued functions such that
\[
{\cal C}_{jk}={\cal R}_{jk}+i{\cal H}_{jk},\;\;
{\cal C}_{j}={\cal R}_{j}+i{\cal H}_{j},\;\;
{\cal C}_{0}={\cal R}_{0}+i{\cal H}_{0}.
\]
We use the notation
\[
{\fr R}(x, D_x)=\sum ^n_{j,k=1}{\cal
R}_{jk}(x)\frac{\partial ^2}{\partial x_j\partial x_k}-
\sum ^n_{j=1}{\cal R}_{j}(x)\frac{\partial}{\partial x_j}-{\cal R}_{0}(x),
\]
\[
{\fr H}(x, D_x)=\sum ^n_{j,k=1}{\cal
H}_{jk}(x)\frac{\partial ^2}{\partial x_j\partial x_k}-
\sum ^n_{j=1}{\cal H}_{j}(x)\frac{\partial}{\partial x_j}-{\cal H}_{0}(x).
\]

Separating the real and imaginary parts of the system ${\fr C}(x, \partial
/\partial x)\bs u=\bs 0$, where $\bs u=\bs v+i\bs w$, we get the following system with real coefficients,
\[
{\fr R}(x, D_x)\bs v-{\fr H}(x, D_x)\bs w=\bs 0,\;\;\;\;
{\fr H}(x, D_x)\bs v+{\fr R}(x, D_x)\bs w=\bs 0,
\]
which, like the original system, is strongly elliptic.

All the assertions below in this subsection are corollaries of the corresponding results in subsection \ref{SELMGR_1.01}.
The next statement is analogous to Theorem \ref{TELMGR_1.07}.

\setcounter{theorem}{9}
\begin{theorem} \label{TELMGR_2.02}
The classical maximum modulus principle
\[
\nl\bs u\nr_{[{\bf C}(\overline \omega)]^m} \leq
\nl\bs u |_{\partial \omega }\nr _{[{\bf C}(\partial  \omega)]^m}
\]
is valid for solutions of the system ${\fr C}(x, D_x)\bs u=\bs 0$ in an arbitrary domain
$\omega \subset \Omega $ with boundary from  the class ${\rm C}^{2, \alpha }$ if and only if:

{\rm (i)} for all $x \in \overline \Omega $ the equalities
\[
{\cal C}_{jk}(x)={\cal C}(x)a_{jk}(x),\;\;1\leq j, k \leq n,
\]
hold, where ${\cal C}$ is a complex $(m\times m)$-matrix-valued function such that
$\Re ({\cal C}(x) \bs \zeta , \bs\zeta ) > 0$ for all $x \in \overline \Omega ,
\bs \zeta \in \bb C^m \backslash \{ \bs 0 \}, ((a_{jk}))$ is a real positive-definite
$(n\times n)$-matrix-valued function in $x \in \overline \Omega $ of order n;

{\rm (ii)}  for all $x \in \Omega $ and any $\bs \xi _j, \bs \zeta \in \bb C^m, j=1,\dots,n$, such that
$\Re (\bs \xi _j, \bs \zeta )=0$ the inequality
\[
\Re \;\Biggl \{\sum ^n_{j,k=1}a_{jk}(x) (\bs\xi _j, \bs\xi _k)+
\sum ^n_{j=1}({\cal C}^{-1}(x)
{\cal C}_j(x)\bs\xi_j, \bs\zeta )+ ({\cal C}^{-1}(x){\cal C}_0(x)
\bs\zeta , \bs\zeta ) \Biggr \} \geq 0
\]
is valid.
\end{theorem}

The following assertion is a consequence of Theorem \ref{TELMGR_2.02}.
\setcounter{theorem}{1}
\begin{corollary} \label{CELMGR_2.03}
The classical maximum modulus principle
\[
\nl\bs u\nr_{[{\bf C}(\overline \omega)]^m} \leq
\nl\bs u |_{\partial \omega }\nr_{[{\bf C}(\partial  \omega)]^m}
\]
holds for solutions of the system
\[
\sum_{j,k=1}^{n} {\cal C}_{jk}(x)\frac{\partial ^2 \bs u}{ \partial x_j \partial x_k}-
\sum_{j=1}^{n} {\cal C}_j(x)\frac{\partial \bs u}{\partial x_j}=\bs 0
\]
in an arbitrary domain $\omega \subset \Omega $ with boundary from the class 
${\rm C}^{2, \alpha }$ if and only if
\[
{\cal C}_{jk}(x)={\cal C}(x)a_{jk}(x),\;\;\;{\cal C}_{j}(x)={\cal C}(x)a_{j}(x),\;\; 1\leq j,\; k\leq n.
\]
Here ${\cal C}(x)$ and $((a_{jk}))$ are the matrix-valued functions 
defined in Theorem $\ref{TELMGR_2.02}$ and
$a_j$ are real scalar functions with the same smoothness as ${\cal C}_j$.
\end{corollary}

\setcounter{theorem}{2}
\begin{remark} \label{R3}  As in Remark \ref{R1}, condition (ii) in Theorem \ref{TELMGR_2.02} 
can be replaced by the following one:
 
{\it for all} $x\in \Omega$ {\it and for any}
$\bs\zeta \in {\bb C}^m\backslash\lbrace \bs 0 \rbrace$ {\it the inequality}
\begin{eqnarray*}
& &|\bs\zeta|^{-2}\sum ^n_{j,k=1}b_{jk}(x)\Re({\cal C}^{-1}(x){\cal C}_j(x)
\bs\zeta , \bs\zeta)({\cal C}^{-1}(x){\cal C}_k(x)\bs\zeta , \bs\zeta ) \\
& &- \sum ^n_{j,k=1}b_{jk}(x)({\cal C}_j^*(x)({\cal C}^*(x))^{-1}\bs\zeta,\;{\cal C}_k^*(x)
({\cal C}^*(x))^{-1}\bs\zeta )\\
& &\\
& &+4\Re({\cal C}^{-1}(x){\cal C}_0(x)\bs\zeta, \bs\zeta )\geq 0
\end{eqnarray*}
{\it is valid, where} $((b_{jk}))$ {\it is the inverse matrix of} $((a_{jk}))$ {\it and} ${\cal C}^*_j(x)$
{\it is the adjoint matrix of} ${\cal C}_j(x)$.
\end{remark}

For the scalar uniformly elliptic equation with complex coefficients of the general form
\begin{equation}\label{SEEQ} 
\sum ^n_{j,k=1}c_{jk}(x)\frac{\partial ^2u}{\partial x_j\partial x_k}-
\sum ^n _{j=1}c_j(x)\frac{\partial u}{\partial x_j}-c_0(x)u=0, 
\end{equation}
Theorem \ref{TELMGR_2.02} and Remark \ref{R3} imply

\setcounter{theorem}{2}
\begin{corollary} \label{CELMGR_2.05}
The classical maximum modulus principle
\[
\nl u\nr_{{\bf C}(\overline \omega) } \leq
\nl u |_{\partial \omega }\nr_{{\bf C}(\partial  \omega)}
\]
is valid for solutions of equation $(\ref{SEEQ})$ in an arbitrary subdomain
$\omega$ of $\Omega$  with boundary $\partial\omega$  of the
class ${\rm C}^{2,\alpha}$ if and only if for all $x\in \overline\Omega$:

{\rm (i)} $\;c_{jk}(x)=c(x)a_{jk}(x),\; 1\leq j,\; k\leq n$,
where $\Re c(x)> 0$ and $((a_{jk}))$  is a real positive-definite $(n\times n)$-matrix-valued function;

{\rm (ii)} the inequality
\[
4\Re \left (\frac{c_0(x)}{c(x)}\right )\geq \sum ^n_{j,k=1}b_{jk}(x)
\Im  \left (\frac{c_j(x)}{c(x)}\right )\Im  \left (\frac{c_k(x)}{c(x)}\right )
\]
holds, where $((b_{jk}))$ is the $(n\times n)$-matrix
inverse of $((a_{jk}))$.
\end{corollary}

\section{Sharp Agmon-Miranda estimates for the gradients of solutions to higher order elliptic equations}\label{S_3}

Everywhere in this section, by smoothness we mean the membership in $\rm C^\infty $.
Suppose that $\Omega$ is a domain in ${\bb R}^n$ with smooth boundary $\partial\Omega$
and a compact closure $\overline\Omega$. We consider the elliptic operator
$$
P(D_x)=\sum _{|\beta|\leq 2\ell}a _{_\beta}D_x^{\beta}
$$
with constant complex coefficients, where $D_x^{\beta}={{\partial ^{|\beta|}}
/ {\partial x_1 ^{\beta _1}\dots\partial x_n ^{\beta _n}}}$, and
$\beta=(\beta_1,$ $\dots,\beta_n)$ is a multi-index of order
$|\beta|=\beta_1+\dots+\beta_n$. By $P_{0}(\xi)$ we denote the
principal homogeneous part of the polynomial $P(\xi)$. For $n=2$ we assume also that all $\xi_2$-roots of the polynomial $P_{0}(\xi)$
for all $\xi_1\in {\bb R} \backslash\lbrace 0\rbrace$.

Let ${\bb R}^{n} _{+}(\bs\nu )=\big \lbrace x\in {\bb R}^{n}: (x, \bs\nu) > 0\big \rbrace$, where $\bs\nu$ is a unit vector and let $K(\bs\nu )$ be the best constant in the Agmon-Miranda inequality 
\begin{equation} \label{Eq_2.1}
\sup _{\overline {{\bb R}^n _{+}(\bs\nu )}} |\nabla _{\ell-1} u|\leq K(\bs\nu )
\sup _{\partial{\bb R}^n _{+}(\bs\nu )} |\nabla _{\ell-1} u|.
\end{equation}
Here
$$
|\nabla _{\ell-1} u|=\left (\sum _{|\beta|=\ell-1} \frac{(\ell-1)!}{\beta !}
|D_x^{\beta} u|^2 \right )^{1/2},
$$
and $u$ is an arbitrary solution of the equation $P_{0}(D_x)u=0$, smooth
in $\overline {{\bb R}^n _{+}(\bs\nu )}$ and such that $u(x)=O(|x|^{\ell-1})$ for large $|x|$.

The following assertion gives a best constant in a weak form (\ref{Eq_1.4}) of the Agmon-Miranda inequality.

\begin{theorem} \label{T_2.1}
For any solution of the equation $P(D_x)u=0$,
smooth on $\overline\Omega$, the inequality
\begin{equation} \label{Eq_2.2}
\max _{\overline\Omega} |\nabla _{\ell-1} u|\leq \left (\sup _{|\bs\nu|=1}K(\bs\nu)\;+\;\varepsilon \right )\max _{\partial\Omega}|
\nabla _{\ell-1} u|+c(\varepsilon)\nl u\nr_{L^1(\Omega)}
\end{equation}
is valid, where $\varepsilon$ is any positive number and $c(\varepsilon)$ is a positive constant independent of $u$.
\end{theorem}

In the next theorem we give we the sharp constant $K$ in the C. Miranda inequality 
\begin{equation} \label{Eq_2.3}
\sup _{\overline {{\bb R}^n _{+}}} |\nabla u|\leq K
\sup _{\partial {\bb R}^n _{+}} |\nabla  u|,
\end{equation}
where $u$ is a solution of the biharmonic equation in ${\bb R}^n _{+}=\{ x=(x_1,\dots, x_n)\in{\bb R}^n: x_n>0 \}$ 
from ${\rm C}^{\infty}(\overline {{\bb R}^n _{+}})$ and $u(x)=O(|x|)$ for large $|x|$.

\begin{theorem} \label{T_2.2}
The sharp constant $K$ in inequality $(\ref{Eq_2.3})$ is given by
$$
K=\frac{2\Gamma (\frac{n}{2})}{\sqrt {\pi }\;\Gamma (\frac{n-1}{2})}
\int ^{\pi /2}_0\big [4+n(n-4)\cos ^2\vartheta \big ]^{1/2}\sin ^{n-2}
\vartheta d\vartheta.
$$
In particular, $K=4/\pi$ for $n=2, K=1/2+2\pi\sqrt{3}/9$ for $n=3$ and $K=2$ for $n=4$.  
\end{theorem}
 
The last assertion was proved in \cite{MK} for the case $n=2$ and
in \cite{KM2} for any $n$. 

\section{Sharp Agmon-Miranda estimates for solutions of the Lam\'e, Stokes and planar deformed state systems} \label{S_4}

Polya's example \cite{PO} demonstrated that the best factor in the inequality between the modulus of the elastic displacement inside the three-dimensional ball and its maximum value on the boundary of that ball exceeds 1. 
A similar inequality (\ref{Eq_1.5A}) with coefficient depending on the domain holds for 
domains with smooth boundary,  and this inequality for solutions of the Lam\'e system is 
called Fichera's maximum principle (see Fichera \cite{FI}). 
This principle is a particular case of the maximum principles for general elliptic 
systems (see Agmon, Douglis and Nirenberg \cite{ADN2}, 
Cannarsa \cite{Can}, Schulze \cite{SHU}, Solonnikov \cite{SOL}, Zhou \cite{CZhou}).

There are works on the Agmon-Miranda type
maximum principle for elliptic systems in domains with singularities at the boundary (Maz'ya and Plamenevski\v{\i} \cite{MaPl}, Albinus \cite{Alb2}, 
Maz'ya and Rossmann \cite{MR3,MR2,MR4} and the bibliography there). 

Estimates for the maximum modulus of velocity vector subject to the nonlinear Navier-Stokes system were obtained
by Solonnikov \cite{SOL4} for smooth domains, by Maz'ya and Rossmann \cite{MR1} for polyhedral domains, by Russo \cite{Russo1} for Lipschitz domains. 
Agmon-Miranda maximum principle as well as existence and uniqueness of solutions to
Stokes system and elastostatics were treated by Maremonti and Russo \cite{MarRu1,MarRu2},
and Tartaglione \cite{ATA}. A survey of maximum principles for the elasticity theory is given by Wheeler \cite{Wheel}.

\subsection{The Lam\'e and Stokes systems}

In the half-space ${\bb R}^n _{+},\; n\geq 2,$ let us consider the Lam\'e system
\begin{equation} \label{(2.2)}
\mu \Delta \bs u+(\lambda +\mu )\hbox {grad}\;\hbox {div}\;\bs u=\bs 0,
\end{equation}
and the Stokes system
\begin{equation} \label{(2.3)}
\nu \Delta \bs u-\hbox {grad}\;p=0,\;\;\;\hbox {div}\;\bs u=\bs 0, 
\end{equation}
with the Dirichlet boundary condition
\begin{equation} \label{(2.4)}
\bs u\big |_{x_{n}=0}=\bs f, 
\end{equation}
where $\lambda $ and $\mu $ are the Lam\'e constants, $\nu$ is the kinematic
coefficient of viscosity, $\bs f\in [{\rm C}_b(\partial {\bb R}^n_{+})]^n$, $\bs u=(u_1,\dots ,u_n)$
is the displacement vector of an elastic medium or the velocity vector of a fluid, and $p$ is the pressure in the fluid.

\smallskip
For the solution $\bs u\in [{\rm C}^2({\bb R}^n _{+})]^n\cap [{\rm C}_{\rm b}(\overline {{\bb R}^n _{+}})]^n$
of problems (\ref{(2.2)}),  (\ref{(2.4)}) and  (\ref{(2.3)}),  (\ref{(2.4)}) we have the representation (see Kupradze, Gegelia, 
Basheleishvili and Burchuladze\cite{KGBB}, Ladyzhenskaya \cite{LA})
$$
\bs u(x)=\int _{\partial {\bb R}^n_{+}}U_{\kappa}\left (\frac{y-x}{|y-x|}\right )
\frac{x_n}{|y-x|^n}\bs f(y')dy',
$$
where $x\in {\bb R}^n _{+},\; y=(y', 0),\; y'=(y_1,\dots,y_{n-1})$.
Here $\kappa =1$ for the Stokes system, $\kappa = (\lambda +\mu )
(\lambda +3\mu )^{-1}$ for the Lam\'e system, and $U_{\kappa}$ is the
$(n\times n)$-matrix-valued function on ${\bb S}^{n-1}=
\big \lbrace x\in {\bb R}^n :\;|x|=1 \big \rbrace$ with the entries
$$
\frac{2}{\omega _n}\Big [(1-\kappa )\delta _{ij}+n\kappa
\frac{(y_{i}-x_{i})(y_{j}-x_{j})}{|y-x|^2}\Big ],
$$
$\omega _n$ being the area of the sphere ${\bb S}^{n-1}$.

Now, we give explicit formulas for sharp constant in inequality
(\ref{Eq_1.6A}) for bounded solutions of Lam\'e and Stokes systems
in a half-space.
\setcounter{theorem}{0}
\begin{theorem} \label{elas}
The sharp constant ${\cal K}({\bb R}^n _{+})$ for the Lam\'e and the Stokes systems in 
$$
|\bs u(x)|\leq {\cal K}({\bb R}^n _{+})\sup \{|\bs u(x')|: x' \in \partial {\bb R}^n _{+} \}
$$
has the form
$$
{\cal K}({\bb R}^n _{+})=\frac{2\Gamma (\frac{n}{2})}{\sqrt {\pi }\Gamma (\frac{n-1}{2})}\int ^{\pi /2}_0\big [(1-\kappa )^2+n\kappa (n\kappa -2\kappa +2)
\cos ^2\vartheta \big ]^{1/2}\sin ^{n-2}\vartheta d\vartheta
$$
and the inequality ${\cal K}({\bb R}^n _{+}) > 1$ holds for $\kappa \neq 0$.

In the case $\kappa=1$, i.e., for an $n$-dimensional Stokes system,
$$
{\cal K}({\bb R}^n _{+})=\frac{2}{\sqrt {\pi }}\frac{\Gamma (\frac{n}{2}+1)}{\Gamma (\frac{n+1}{2})}.
$$
\end{theorem}

\setcounter{theorem}{0}
Below we give consequences of Theorem \ref{elas} for the cases
$n=2$ and $n=3$, respectively.
\begin{corollary}
The equality
\begin{eqnarray*}
{\cal K}({\bb R}^2 _{+})&=&\frac{2}{\pi }(1+\kappa )E\left (\frac{2\sqrt {\kappa }}{1+\kappa }\right )\\
&=&1+\frac{1}{2^2}\kappa ^2+\frac{1}{2^2 4^2}
\kappa ^4+\dots + \left [\frac{(2m-3)!!}{2^mm!}\right ]^2\kappa ^{2m}+\dots
\end{eqnarray*}
is valid, where $E$ is the complete elliptic integral of
the second kind. In particular, ${\cal K}({\bb R}^2 _{+})=4/\pi$
for $\kappa=1$.
\end{corollary}

\begin{corollary}
The equality
$$
{\cal K}({\bb R}^3 _{+})=\frac{1}{2}\left [1+2\kappa + \frac{(1-\kappa )^2}{\sqrt {3\kappa (\kappa +2)}}\log \frac{1+2\kappa + \sqrt {3\kappa (\kappa +2)}}{
1-\kappa }\right ]
$$
is valid. In particular, ${\cal K}({\bb R}^3 _{+})=3/2$ for $\kappa=1$.
\end{corollary}

\subsection{Planar deformed state}

Let $\sigma _{11}, \sigma _{12}$ and $\sigma _{22}$ be the components of
the stress tensor in the half-plane ${\bb R}^2 _{+}$. Consider the system of
equations in ${\bb R}^2 _{+}$ for the stresses in a planar deformed state
(see, for example, Muskhelishvili \cite{Musk}):
$$
\partial\sigma _{11}/\partial x_1+\partial\sigma _{12}/\partial x_2=0,
$$
$$
\partial\sigma _{12}/\partial x_1+\partial\sigma _{22}/\partial x_2=0,
$$
$$
\Delta (\sigma _{11}+\sigma _{22})=0,
$$
with the boundary conditions
$$
\sigma _{12}(x_1, 0)=p_1(x_1),\;\;\sigma _{22}(x_1, 0)=p_2(x_1),
$$
where $p_1$ and $p_2$ are continuous and bounded functions on $\partial
{\bb R}^2_{+}$.

\begin{theorem} The sharp constant in the inequality
$$
||(\sigma _{12}^2+\sigma _{22}^2)^{1/2}|| _{C(\overline {{\bb R}^2_{+}})} \leq
{\cal K}||(\sigma _{12}^2+\sigma _{22}^2)^{1/2}|| _{C(\partial {\bb R}^2_{+})}
$$
is equal to $4/\pi$.
\end{theorem}

\medskip
{\bf Acknowledgement.} The publication has been prepared with the support of the "RUDN University Program 5-100".



\begin{thebibliography}{99}

\bibitem{AG0} S. Agmon, \textit{ Multiple layer potentials and the Dirichlet 
problem for higher order elliptic equations in the plane, I}, 
Comm. Pure Appl. Math., \textbf{10} (1957), 179--239.

\bibitem{AG} S. Agmon, \textit{Maximum theorems for solutions of higher order elliptic equations}, Bull. Amer. Math. Soc., \textbf{66} (1960), 77--80.

\bibitem{ADN1} S. Agmon, A. Douglis, and L. Nirenberg,
\textit{Estimates near the boundary for solutions of elliptic partial
differential equations satisfying general boundary conditions, I.},
Comm. Pure Appl. Math., \textbf{12} (1959), 623--727.

\bibitem{ADN2} S. Agmon, A. Douglis, and L. Nirenberg,
\textit{Estimates near the boundary for solutions of elliptic partial
differential equations satisfying general boundary conditions, II.},
Comm. Pure Appl. Math., \textbf{17} (1964), 35--92.

\bibitem{Alb2} G. Albinus, \textit{Estimates of Miranda-Agmon type in plane domains with corners 
and existence theorem}, in: Constructive Function Theory '81,  Sofia, 197--204.

\bibitem{BI} A.V. Bitsadze, \textit{On elliptic systems of second order
partial differential equations}, Dokl. Akad. Nauk SSSR, \textbf{112} (1957),
983--986 (Russian).

\bibitem{BITS} A.V. Bitsadze,  \textit{Boundary value problems for second
order elliptic equations}, Moskva, Nauka, 1966;
English translation: North-Holland Publ. Company, Amsterdam, 1968.

\bibitem{Can} P. Cannarsa, \textit{On a maximum principle for elliptic systems with constant coefficients}, Rend. Sem. Mat. Univ. Padova, \textbf{64} (1981), 77--84.

\bibitem{FI} G. Fichera, \textit{Il teorema del massimo modulo per l'equazione
dell'elastostatica tridimensionale}, Arch. Rat. Mech. Anal., {\bf 7} (1961),
373--387.

\bibitem{FiMi} D.G. de Figueiredo and E. Mitidieri, \textit{Maximum principles for linear
elliptic systems}, Rend. Istit. Mat. Univ. Trieste,  {\bf 22}:1-2 (1992),
36--66.

\bibitem{HP} G. Hile and M.H. Protter, \textit{Maximum principles
for a class of first-order elliptic systems}, J. Diff. Equations,
{\bf 24} (1977), 136--151.

\bibitem{Hong} C.W. Hong, \textit{Necessary and sufficient condition for a class of 
generalized maximum principles}, Acta Math. Sinica, {\bf 26}:3 (1983), 307--321.
(chinese).

\bibitem{KamKh1} L.I. Kamynin and B.N. Khimchenko, \textit{On the weak maximum principle for 
a second order elliptic system}, Izv. Ross. Akad. Nauk, Ser. Mat., {\bf 59}:5 (1995),
73--84 (Russian); English. transl. in Izv. Math., {\bf 59}:5 (1995), 949--961.

\bibitem{KamKh2} L.I. Kamynin and B.N. Khimchenko, \textit{Necessary and sufficient conditions
for satisfying the weak extremum principle for second order elliptic systems}, 
Sibirsk. Mat. Zh., {\bf 37}:6 (1996), 1314--1334 (Russian); English. transl.
in Siberian Math. J., {\bf 37}:6  (1996), 1153--1170.

\bibitem{KM2} G.I. Kresin and V.G. Maz'ya, \textit{On the exact
constant in the inequality of Miranda-Agmon type for solutions of
elliptic equations}, Izvestia Vys. Uch. Zaved., Matem. Ser.,
{\bf 5} (1988), 41-50 (Russian); English transl.: Izvestia VUZ. Matematika, 
{\bf 32}:5 (1988), 49--59.

\bibitem{KM3} G.I. Kresin and V.G. Maz'ya, \textit{Criteria for
validity of the maximum modulus principle for solutions of linear
strongly elliptic systems}, Potential Analysis, {\bf 2} (1993), 73--99.

\bibitem{KM} G. Kresin and V. Maz'ya, \textit{Maximum Principles and Sharp 
Constants for Solutions of Elliptic and Parabolic Systems}, Math. Surveys and Monographs, {\bf 183}, Amer. Math. Soc., Providence, 
Rhode Island, 2012.

\bibitem{KGBB} V.D. Kupradze, T.G. Gegelia, M.O. Basheleishvili, and T.V. Burchuladze, 
\textit{Three-dimensional Problems of the Mathematical
Theory of Elasticity and Thermoelasticity}, North-Holland, Amsterdam, 1979.

\bibitem{LA} O.A. Ladyzhenskaya, \textit{The Mathematical Theory
of Viscous Incompressible Flow}, Gordon and Breach, New York, 1963.

\bibitem{LeS} S. Lenhart and P. Schaefer \textit{On Comparison results for classical and
viscosity solutions in elliptic systems}, in: ``Evolution Equations'', G. Ferreyra,
G. R. Goldstein, F. Neubrander ed., Marcel Dekker, 1995, pp. 269--275.

\bibitem {LO1} Ya.B. Lopatinski\v{\i}, \textit{On a method of reducing boundary value
problems for systems of differential equations of elliptic type to regular integral
equations}, Ukrain. Mat. \v{Z}urnal, {\bf 5}:2 (1953), 123--151 (Russian).

\bibitem{LoMo} J. L\'opez-G\'omez and M. Molina-Meyer, \textit{The maximum principle for 
cooperative weakly coupled elliptic systems and some applications}, Diff. Int.
Eq., {\bf 7}:2 (1994), 383--398.

\bibitem{MarRu1} P. Maremonti and R. Russo, \textit{On the maximum modulus theorem
for the Stokes system}, Ann. Scuola Norm. Sup. Pisa (IV), {\bf 30} (1994), 630--643.

\bibitem{MarRu2} P. Maremonti and R. Russo, \textit{On existence and uniqueness
of classical solutions of the stationary Navier-Stokes equations and to the
traction problem of linear elastostatics}, Quaderni di Matematica. 
Classical Problems in Mechanics, {\bf 1} (1997), 171--251.

\bibitem{MK} V.G. Maz'ya and G.I. Kresin, \textit{On the maximum
principle for strongly elliptic and parabolic second order systems
with constant coefficients}, Mat. Sb., {\bf 125(167)} (1984), 458--480 (Russian);
English transl.: Math. USSR Sb., {\bf 53} (1986) ,457--479.

\bibitem {MaPl} V. Maz'ya and B.A. Plamenevski\v{\i}, \textit{Estimates in $L_p$ and H\"older classes
and the Miranda-Agmon maximum principle for solutions of elliptic boundary value problems in domains with singular points on the boundary}, Math. Nachr., {\bf 81} (1978), 25--82, Engl. transl. in: Amer. Math. Soc. Transl., {\bf 123} (1984), 1--56.

\bibitem{MR3} V. Maz'ya and J. Rossmann, \textit{On the Agmon-Miranda maximum principle 
for solutions of elliptic equations in polyhedral and polygonal domains}, 
Ann. Global Anal. Geom., {\bf 9} (1991), 253--303.

\bibitem{MR2} V. Maz'ya and J. Rossmann, \textit{On the Agmon-Miranda maximum principle 
for solutions of strongly elliptic equations in domains of ${\mathbb R}^n$ 
with conical points}, Ann. Global Anal. Geom., {\bf 10} (1992), 125--150.

\bibitem{MR1} V. Maz'ya and J. Rossmann, \textit{A maximum modulus estimate for solutions
of the Navier-Stokes system in domains of polyhedral type}, Math. Nachr., {\bf 282}:3 (2009), 459--469.

\bibitem{MR4} V. Maz'ya and J. Rossmann, \textit{Elliptic Equations in Polyhedral Domains}, 
Mathematical Surveys and Monographs, {\bf 162}, Amer. Math. Soc., RI, 2010.

\bibitem{MIR1} C. Miranda, \textit{Formule di maggiorazione e
teorema di esistenza per le funzioni biarmoniche di due variabili},
Giorn. Mat. Bataglini, {\bf 78} (1948-1949), 97--118.

\bibitem{MIR2} C. Miranda, \textit{Teorema del massimo modulo e
teorema di esistenza e di unicit\`a per il problema di Dirichlet
relativo alle equazioni ellitiche in due variabili}, Ann. Mat. Pura
Appl., ser. 4, {\bf 46} (1958), 265--311.

\bibitem{MIR3} C. Miranda, \textit{Sul teorema del massimo modulo per
una classe di sistemi ellitici di equazioni del secondo ordine e per
le equazioni a coefficienti complessi}, Istit. Lombardo Accad. Sci.
Lett. Rend. A, {\bf 104} (1970), 736--745.

\bibitem{MiSw} E. Mitidieri and G. Sweers, \textit{Weakly coupled elliptic systems and
positivity}, Math. Nachr., {\bf 173}, (1995), 259--286.
 
\bibitem{Musk} N.I. Muskhelishvili, \textit{Some Basic Problems of the Mathematical
Theory of Elasticity}, Noordhoff, Groningen-Holland, 1953.

\bibitem{PI} B. Pini, \textit{Sui sistemi equazioni lineari a derivate
parziali del secondo ordine dei tipi ellittico e parabolico}, Rend.
del Sem. Mat. della Univ. di Padova, {\bf 22} (1953), 265-280.
 
\bibitem{PO} G. Polya, \textit{Liegt di Stelle der gr\"o$\beta$ten
Beanspruchungen an der Oberfl\"ache?}, Z. Angew. Math. Mech.,
{\bf 10} (1930), 353--360. 
 
\bibitem{PW} M.H. Protter and H.F. Weinberger, \textit{Maximum
Principles in Differential Equations}, Prentice-Hall, Inc., Englewood Cliffs, N.J., 1967;
Springer-Verlag, New York Inc., 1984.

\bibitem{PR}  M.H. Protter, \textit{Maximum principles}, in ``Maximum principles 
and eigenvalue problems in partial differential equations'',  P.W. Schaefer ed., Pittman
Research Notes, Math. ser., {\bf 175} (1988), pp. 1--14.

\bibitem{Rus} I.A. Rus, \textit{Un principe du maximum pour les solutions d'un
systeme fortement elliptique}, Glasnik Mat., {\bf 4(24)}:1 (1969), 75-77.

\bibitem{Russo1} R. Russo, \textit{On the maximum modulus theorem
for the steady-state Navier-Stokes equations in Lipschitz bounded domains}, 
Appl. Anal., {\bf 90}:1-2 (2011), 193--200.

\bibitem{SAB2} K.B. Sabitov, \textit{Maximum modulus principle for some classes of second order
elliptic and hyperbolic systems}, Differentsial'nye uravnenia, {\bf 27}:2 (1991), 
272--278 (Russian); English. transl. in Differential Equations, {\bf 27}:2, (1991),
194--199.

\bibitem{SHU} B.-W. Schulze, \textit{On a priori estimates in uniform norms for strongly
elliptic systems}, Sibirsk. Mat. Zh., {\bf 16}:2 (1975), 384--394 (Russian); English. transl.
in Siberian Math. J., {\bf 16}:2 (1975), 297--305.

\bibitem{SHA} Z.Ya. Shapiro, \textit{The first boundary value problem for an
elliptic system of differential equations}, Mat. Sb., {\bf 28(70)}:1 (1951),
55-78 (Russian).

\bibitem{SirB} B. Sirakov, \textit{Some estimates and maximum principles for weakly coupled
systems of elliptic PDE }, Nonlinear Anal., {\bf 70}:8 (2009), 3039--3046.

\bibitem{SOL} V.A. Solonnikov, \textit{On general boundary value problems
for systems elliptic in the Douglis-Nirenberg sense}, Izv. Akad. Nauk SSSR, ser.
Mat., {\bf 28}:3 (1964), 665--706 ; English transl. Amer. Math. Soc. Transl. (2), 
{\bf 56} (1966), 193--232.

\bibitem{SOL4} V.A. Solonnikov, \textit{On a maximum modulus estimate of the solution
of stationary problem for the Navier-Stokes equations}, Zap. Nauchn. Sem. S.-Petersburg
Otdel. Mat. Inst. Steklov (POMI), {\bf 249} (1997), 
294--302; Transl. in J. Math. Sci. (New York), {\bf 101}:5 (2000), 3563--3569.
 
\bibitem{STY1} T. Stys, \textit{Aprioristic estimations of solutions of a 
certain elliptic system of partial differential second order
equations}, Bull. Acad. Pol. Sc., {\bf 13} (1965), 639--640.

\bibitem{SZE} P. Szeptycki, \textit{Existence theorem for the first boundary
value problem for a quasilinear elliptic system}, Bulletin
Acad. Polon. des Sciences, {\bf 7} (1959), 419--424.

\bibitem{ATA} A. Tartaglione, \textit{On existence, uniqueness and the maximum modulus 
theorem in plane linear elastostatics for exterior domains}, Ann. Univ. 
Ferrara - Sez. VII - Sc. Mat., {\bf 47} (2001), 89--106.

\bibitem{WAS} J. Wasowski, \textit{Maximum principles for a certain strongly
elliptic system of linear equations of second order},
Bull. Acad. Polon. Sci., {\bf 18} (1970), 741--745.

\bibitem{Wheel} L.T. Wheeler, \textit{Maximum principles in classical elasticity},
in: ``Mathematical problems in elasticity'', R. Russo ed., 
Singapore World Sc., (1996), pp. 157--185.

\bibitem{CZhou} C. Zhou, Maximum principles for elliptic systems, 
{\it J. Math. Anal. and Appl.}, {\bf 159} (1991), 418--439.

\end{thebibliography}
\end{document}